\documentclass[12pt,twoside]{article}
\usepackage[latin1]{inputenc}
\usepackage{amsmath}
\usepackage{amssymb,amsfonts}
\usepackage{graphicx}
\usepackage{times,amssymb,amscd}
\usepackage[USenglish]{babel}
\usepackage{verbatim}
\usepackage{enumerate}

\newcommand{\bH}{\mathbf{H}}

\newcommand{\bL}{\mathbf{L}}

\newcommand{\bR}{\mathbf{R}}

\newcommand{\bS}{\mathbf{S}}

\newcommand{\bt}{\mathbf{t}}

\newcommand{\cD}{\mathcal{D}}

\newcommand{\cP}{\mathcal{P}}

\newcommand{\cC}{\mathcal{C}}
\newcommand{\cB}{\mathcal{B}}

\newcommand{\SXR}{\bS^2\!\times\!\bR}
\newcommand{\HXR}{\bH^2\!\times\!\bR}
\newcommand{\SLR}{\widetilde{\bS\bL_2\bR}}
\newcommand{\NIL}{\mathbf{Nil}}
\newcommand{\SOL}{\mathbf{Sol}}

\begin{document}
\pagestyle{myheadings}
\markboth{\centerline{J.~Szirmai}}
{Horoball packings related to hyperbolic 24 cell}
\title
{Horoball packings \\ related to hyperbolic $24$ cell
\footnote{Mathematics Subject Classification 2010: 52C17, 52C22, 52B15. \newline
Key words and phrases: Hyperbolic geometry, horoball packings, polyhedral density function, optimal density.
}}

\author{\medbreak \medbreak {\normalsize{}} \\
\normalsize Jen\H{o}  Szirmai \\
\normalsize  Budapest University of Technology and Economics\\
\normalsize Institute of Mathematics, Department of Geometry \\
\normalsize H-1521 Budapest, Hungary \\
\normalsize Email: szirmai@math.bme.hu }
\date{\normalsize (\today)}
\maketitle
\begin{abstract}

In this paper we study the horoball packings related to the hyperbolic 24 cell in the extended hyperbolic space $\overline{\mathbf{H}}^4$
where we allow {\it horoballs in different types}
centered at the various vertices of the 24 cell.

We determine, introducing the notion of the generalized polyhedral density function, the locally densest horoball packing arrangement
and its density with respect to the above regular tiling. The maximal density is $\approx 0.71645$ which is equal to the known
greatest ball packing density in hyperbolic 4-space given in \cite{KSz14}.
\end{abstract}


\newtheorem{theorem}{Theorem}[section]
\newtheorem{corollary}[theorem]{Corollary}
\newtheorem{lemma}[theorem]{Lemma}
\newtheorem{exmple}[theorem]{Example}
\newtheorem{defn}[theorem]{Definition}
\newtheorem{rmrk}[theorem]{Remark}
\newtheorem{proposition}[theorem]{Proposition}
\newenvironment{definition}{\begin{defn}\normalfont}{\end{defn}}
\newenvironment{remark}{\begin{rmrk}\normalfont}{\end{rmrk}}
\newenvironment{example}{\begin{exmple}\normalfont}{\end{exmple}}
\newenvironment{acknowledgement}{Acknowledgement}




\section{Introduction}

We consider horospheres and their bodies, the horoballs. A horoball packing $\cB$ of $\overline{\mathbf{H}}^n$ is an
arrangement of non-overlapping horoballs ${B}$ in $\overline{\mathbf{H}}^n$.

The definition of packing density is critical in hyperbolic space as shown by B\"or\"oczky \cite{B78}. For standard examples also see
\cite{R06}.
The most widely accepted notion of packing density considers the local densities of
balls with respect to their Dirichlet-Voronoi cells (cf. \cite{B78} and \cite{K98}).
In order to consider horoball packings in $\overline{\mathbb{H}}^n$ we use an extended notion of such local density.

Let $B$ be a horoball in packing $\cB$, and $P \in \overline{\mathbb{H}}^n$ be an arbitrary point.
Define $d(P,B)$ to be the perpendicular distance from point $P$ to the horosphere $S = \partial B$, where $d(P,B)$
is taken to be negative when $P \in B$. The Dirichlet--Voronoi cell $\cD(B,\cB)$ of a horoball $B$ of packing $\cB$ is defined as the convex body
\begin{equation}
\cD(B,\cB) = \{ P \in \mathbb{H}^n | d(P,B) \le d(P,B'), ~ \forall B' \in \cB \}. \notag
\end{equation}
Both $B$ and $\cD$ are of infinite volume, so the usual notion of local density is
modified as follows. Let $Q \in \partial{\mathbb{H}}^n$ denote the ideal center of $B$ at infinity, and take its boundary $S$ to be the one-point compactification of Euclidean $(n - 1)$-space.
Let $B_C^{n-1}(r) \subset S$ be an $(n-1)$-ball with center $C \in S \setminus \{Q\}$.
Then $Q \in \partial {\mathbb{H}^n}$ and $B_C^{n-1}(r)$ determine a convex cone
$\cC^n(r) = cone_Q(B_C^{n-1}(r)) \in \overline{\mathbb{H}}^n$ with
apex $Q$ consisting of all hyperbolic geodesics passing through $B_C^{n-1}(r)$ with limit point $Q$. The local density $\delta_n(B, \cB)$ of $B$ to $\cD$ is defined as
\begin{equation}
\delta_n(\cB, B) =\varlimsup\limits_{r \rightarrow \infty} \frac{vol(B \cap \cC^n(r))} {vol(\cD \cap \cC^n(r))}. \notag
\end{equation}
This limit is independent of the choice of center $C$ for $B^{n-1}_C(r)$.

For periodic ball or horoball packings the local density defined above can be extended to the entire hyperbolic space.
This local density is related to the simplicial density function that was generalized in \cite{Sz11-2} and \cite{Sz13}.
In this paper we will use the generalization of this definition of packing density.

{\it In \cite{Sz11-2}  we have refined the notion of the ,,congruent" horoballs in a horoball packing to the horoballs of the ,,same type"
because the horoballs are always congruent in the hyperbolic space $\overline{\mathbf{H}}^n$, in general.

Two horoballs in a horoball packing are in the ,,same type", or ,,equipacked", if and only if the
local densities of the horoballs to the corresponding cell (e.g. D-V cell; or ideal regular polytop, later on) are equal.}

{\it \bf If we assume that the ,,horoballs belong to the same type"}, then by analytical continuation,
the well known simplicial density function on $\overline{\mathbf{H}}^n$ can be extended from $n$-balls of radius $r$ to the case $r = \infty$,
too. Namely, in this case consider $n + 1$ horoballs which are mutually tangent and let $B$ be one of them. The convex hull of their
base points at infinity will be a totally asymptotic or ideal regular simplex $T_{reg}^{\infty} \in \overline{\mathbf{H}}^n$ of finite volume.
Hence, in this case it is legitimated to write
\begin{equation}
d_n(\infty) = (n + 1)\frac{vol(B \cap T_{reg}^\infty)}{vol(T_{reg}^\infty)}. \notag
\end{equation}
Then for a horoball packing $\cB$, there is an analogue of ball packing, namely (cf. \cite{B78}, Theorem
4)
\begin{equation}
\delta_n(\cB, B) \le d_n(\infty),~ \forall B \in \cB. \notag
\end{equation}
\begin{rmrk}
The upper bound $d_n(\infty)$ $(n=2,3)$ is attained for a regular horoball packing, that is, a
packing by horoballs which are inscribed in the cells of a regular honeycomb of $\overline{\mathbf{H}}^n$. For
dimensions $n = 2$, there is only one such packing. It belongs to the regular tessellation $\{\infty, 3 \}$ . Its dual
$\{3,\infty\}$ is the regular tessellation by ideal triangles all of whose vertices are surrounded
by infinitely many triangles. This packing has in-circle density $d_2(\infty)=\frac{3}{\pi} \approx 0.95493 $.

In $\overline{\mathbf{H}}^3$ there is exactly one horoball packing with horoballs in same type whose Dirichlet--Voronoi cells give rise to a
regular honeycomb described by the Schl\"afli symbol $\{6,~3,~3\}$ . Its
dual $\{3,3,6\}$ consists of ideal regular simplices $T_{reg}^\infty$  with dihedral angle $\frac{\pi}{3}$ building up a 6-cycle around each edge
of the tessellation. The density of this packing is $\delta_3^\infty\approx 0.85328$
\end{rmrk}

{\bf If horoballs of different types at the various ideal vertices are allowed} i.e the horoballs are differently packed, then we generalized
the notion of the simplicial density function \cite{Sz11-2}.
In \cite{KSz} we proved that the optimal ball packing arrangement in $\mathbb{H}^3$ mentioned above is not unique.
We gave several new examples of horoball packing arrangements based on totally asymptotic Coxeter tilings that yield the
B\"or\"oczky--Florian upper bound \cite{BF64}.

Furthermore, in \cite{Sz11-2}, \cite{Sz13} we found that
by admitting horoballs of different types at each vertex of a totally asymptotic simplex and generalizing
the simplicial density function to $\overline{\mathbf{H}}^n$ for $(n \ge 2)$,  the B\"or\"oczky-type density
upper bound is no longer valid for the fully asymptotic simplices for $n \ge 3$.
For example, in $\overline{\mathbb{H}}^4$ the locally optimal packing density was found to be 
$0.77038\dots$ which is higher than the B\"or\"oczky-type density upper bound $0.73046\dots$.
However these ball packing configurations are only locally optimal and cannot be extended to the entirety of the
hyperbolic spaces $\mathbb{H}^n$.

In \cite{KSz14} we have continued our investigations
on ball packings in hyperbolic 4-space. Using horoball packings, allowing horoballs of different types,
we find seven counterexamples with density $\approx 0.71645$ (which are realized by allowing up to three horoball types)
to one of L. Fejes-T\'oth's conjectures.

Several extremal properties relate to the regular hyperbolic 24-cell and the corresponding
Coxeter honeycomb concerning the right angled polytops and hyperbolic 4-manifolds.

A. Kolpakov in \cite{Ko14} has shown that the hyperbolic
24-cell has minimal volume and minimal facet number among all ideal right-angled polytopes in $\overline{\mathbf{H}}^4$.

J.~G. Ratcliffe and S.~T. Tschantz in \cite{RS} have constructed complete, open, hyperbolic 4-manifolds of smallest volume by
gluing together the sides of a regular ideal 24-cell
in hyperbolic 4-space. They also showed that the volume spectrum of hyperbolic 4-manifolds is the set
of all positive integral multiples of $4\pi^2/3$.

L. Slavich has constructed in \cite{S14}, using the hyperbolic 24-cell, two new examples of non-orientable, noncompact,
hyperbolic 4-manifolds. The first has minimal volume $V_m=4\pi^2/3$
and two cusps. This example has the lowest number of cusps among known
minimal volume hyperbolic 4-manifolds. The second has volume $2\cdot V_m$ and
one cusp. It has lowest volume among known one-cusped hyperbolic 4-manifolds.

{\it In this paper we study a new extremal property of the hyperbolc regular 24-cell and the corresponding
regular $4$-dimensional honeycomb described by the Schl\"afli symbol $\{3,~4,~3,~4\}$
relating to horoball packings.

We determine, introducing the notion of the generalized polyhedral density function, the locally densest horoball packing arrangements
and their densities with respect to the above 4-dimensional regular tiling. The maximal density is $\approx 0.71645$ which is equal to the known
greatest ball packing density in hyperbolic 4-space given in \cite{KSz14}.}

\section{Formulas in the projective model}

We use the projective model in Lorentzian $(n+1)$-space
$\mathbb{E}^{1,n}$ of signature $(1,n)$, i.e.~$\mathbb{E}^{1,n}$ is
the real vector space $\mathbf{V}^{n+1}$ equipped with the bilinear
form of signature $(1,n)$
\begin{equation}
\langle ~ \mathbf{x},~\mathbf{y} \rangle = -x^0y^0+x^1y^1+ \dots + x^n y^n \label{bilinear_form} \tag{2.1}
\end{equation}
where the non-zero real vectors
$\mathbf{x}=(x^0,x^1,\dots,x^n)\in\mathbf{V}^{n+1}$
and $ \mathbf{y}=(y^0,y^1,\dots,$ $y^n)\in\mathbf{V}^{n+1}$ represent points in projective space
$\mathcal{P}^n(\mathbb{R})$. $\mathbb{H}^n$ is represented as the
interior of the absolute quadratic form
\begin{equation}
Q=\{[\mathbf{x}]\in\mathcal{P}^n | \langle ~ \mathbf{x},~\mathbf{x} \rangle =0 \}=\partial \mathbb{H}^n  \tag{2.2}
\end{equation}
in real projective space $\mathcal{P}^n(\mathbf{V}^{n+1},\mbox{\boldmath$V$}\!_{n+1})$. All proper interior points $\mathbf{x} \in \mathbb{H}^n$ are characterized by
$\langle ~ \mathbf{x},~\mathbf{x} \rangle < 0$.

The boundary points $\partial \mathbb{H}^n $ in
$\mathcal{P}^n$ represent the absolute points at infinity of $\mathbb{H}^n$.
Points $\mathbf{y}$ satisfying $\langle ~ \mathbf{y},~\mathbf{y} \rangle >
0$ lie outside $\partial \mathbb{H}^n $ and are called the outer points
of $\mathbb{H}^n$. Take $P([\mathbf{x}]) \in \mathcal{P}^n$, point
$[\mathbf{y}] \in \mathcal{P}^n$ is said to be conjugate to
$[\mathbf{x}]$ relative to $Q$ when $\langle ~
\mathbf{x},~\mathbf{y} \rangle =0$. The set of all points conjugate
to $P([\mathbf{x}])$ form a projective (polar) hyperplane
\begin{equation}
pol(P):=\{[\mathbf{y}]\in\mathcal{P}^n | \langle ~ \mathbf{x},~\mathbf{y} \rangle =0 \}. \tag{2.3}
\end{equation}
Hence the bilinear form $Q$ in (\ref{bilinear_form}) induces a bijection
or linear polarity $\mathbf{V}^{n+1} \rightarrow
\mbox{\boldmath$V$}\!_{n+1}$
between the points of $\mathcal{P}^n$
and its hyperplanes.
Point $X [\bold{x}]$ and hyperplane $\alpha
[\mbox{\boldmath$a$}]$ are incident if the value of
linear form $\mbox{\boldmath$a$}$ evaluated on vector $\bold{x}$ is
 zero, i.e. $\bold{x}\mbox{\boldmath$a$}=0$ where $\mathbf{x} \in \
\mathbf{V}^{n+1} \setminus \{\mathbf{0}\}$, and $\ \mbox{\boldmath$a$} \in
\mbox{\boldmath$V$}_{n
+1} \setminus \{\mbox{\boldmath$0$}\}$.
Similarly, lines in $\mathcal{P}^n$ are characterized by
2-subspaces of $\mathbf{V}^{n+1}$ or $(n-1)$-spaces of $\
\mbox{\boldmath$V$}\!_{n+1}$ \cite{Mol97}.

Let $P \subset \mathbb{H}^n$ denote a polyhedron bounded by
a finite set of hyperplanes $H^i$ with unit normal vectors
$\mbox{\boldmath$b$}^i \in \mbox{\boldmath$V$}\!_{n+1}$ directed
 towards the interior of $P$:
\begin{equation}
H^i:=\{\mathbf{x} \in \mathbb{H}^d | \langle ~ \mathbf{x},~\mbox{\boldmath$b$}^i \rangle =0 \} \ \ \text{with} \ \
\langle \mbox{\boldmath$b$}^i,\mbox{\boldmath$b$}^i \rangle = 1. \tag{2.4}
\end{equation}
In this paper $P$ is assumed to be an acute-angled polyhedron
with proper or ideal vertices.
The Grammian matrix $G(P):=( \langle \mbox{\boldmath$b$}^i,
\mbox{\boldmath$b$}^j \rangle )_{i,j} ~ {i,j \in \{ 0,1,2 \dots n \} }$  is an
indecomposable symmetric matrix of signature $(1,n)$ with entries
$\langle \mbox{\boldmath$b$}^i,\mbox{\boldmath$b$}^i \rangle = 1$
and $\langle \mbox{\boldmath$b$}^i,\mbox{\boldmath$b$}^j \rangle
\leq 0$ for $i \ne j$ where
$$
\langle \mbox{\boldmath$b$}^i,\mbox{\boldmath$b$}^j \rangle =
\left\{
\begin{aligned}
&0 & &\text{if}~H^i \perp H^j,\\
&-\cos{\alpha^{ij}} & &\text{if}~H^i,H^j ~ \text{intersect \ along an edge of $P$ \ at \ angle} \ \alpha^{ij}, \\
&-1 & &\text{if}~\ H^i,H^j ~ \text{are parallel in the hyperbolic sense}, \\
&-\cosh{l^{ij}} & &\text{if}~H^i,H^j ~ \text{admit a common perpendicular of length} \ l^{ij}.
\end{aligned}
\right.
$$
This is visualized using the weighted graph or scheme of the polytope $\sum(P)$. The graph nodes correspond
to the hyperplanes $H^i$ and are connected if $H^i$ and $H^j$ not perpendicular ($i \neq j$).
If they are connected we write the positive weight $k$ where  $\alpha_{ij} = \pi / k$ on the edge, and
unlabeled edges denote an angle of $\pi/3$.

In this paper we set the sectional curvature of $\mathbb{H}^n$,
$K=-k^2$, to be $k=1$. The distance $d$ of two proper points
$[\mathbf{x}]$ and $[\mathbf{y}]$ is calculated by the formula
\begin{equation}
\cosh{{d}}=\frac{-\langle ~ \mathbf{x},~\mathbf{y} \rangle }{\sqrt{\langle ~ \mathbf{x},~\mathbf{x} \rangle
\langle ~ \mathbf{y},~\mathbf{y} \rangle }}. \tag{2.5}
\end{equation}
The perpendicular foot $Y[\mathbf{y}]$ of point $X[\mathbf{x}]$ dropped onto plane $[\mbox{\boldmath$u$}]$ is given by
\begin{equation}
\mathbf{y} = \mathbf{x} - \frac{\langle \mathbf{x}, \mathbf{u} \rangle}{\langle \mathbf{u}, \mathbf{u} \rangle} \mathbf{u}, \tag{2.6}
\end{equation}
where $[\mathbf{u}]$ is the pole of the plane $[\mbox{\boldmath$u$}]$.
\medbreak
A horosphere in $\mathbb{H}^n$ ($n \ge 2)$ is a
hyperbolic $n$-sphere with infinite radius centered
at an ideal point on $\partial \mathbb{H}^n$. Equivalently, a horosphere is an $(n-1)$-surface orthogonal to
the set of parallel straight lines passing through a point of the absolute quadratic surface.
A horoball is a horosphere together with its interior.

We consider the usual Beltrami-Cayley-Klein ball model of $\mathbb{H}^n$
centered at $O(1,0,0,$ $\dots, 0)$ with a given vector basis
$\bold{a}_i \ (i=0,1,2,\dots, n)$ and set an
arbitrary point at infinity to lie at $T_0=(1,0,\dots, 0,1)$.
The equation of a horosphere with center
$T_0=(1,0,\dots,1)$ passing through point $S=(1,0,\dots,s)$ is derived from the equation of the
the absolute sphere $-x^0 x^0 +x^1 x^1+x^2 x^2+\dots + x^n x^n = 0$, and the plane $x^0-x^n=0$ tangent to the absolute sphere at $T_0$.
The general equation of the horosphere is in projective coordinates ($s \neq \pm1$):
\begin{align}
(s-1)\left(-x^0 x^0 +\sum_{i=1}^n (x^i)^2\right)-(1+s){(x^0-x^n)}^2 & =0, \tag{2.7}
\end{align}
and in cartesian coordinates setting $h_i=\frac{x^i}{x^0}$ it becomes
\begin{equation}
\label{eqn:horosphere1}
\frac{2 \left(\sum_{i=1}^n h_i^2 \right)}{1-s}+\frac{4 \left(h_d-\frac{s+1}{2}\right)^2}{(1-s)^2}=1. \tag{2.8}
\end{equation}

In $n$-dimensional hyperbolic space any two horoballs are congruent in the classical sense.
However, it is often useful to distinguish between certain horoballs of a packing.
We use the notion of horoball type with respect to the packing as introduced in \cite{Sz11-2}.

{\it Two horoballs of a horoball packing are said to be of the same type or equi\-packed if
and only if their local packing densities with respect to a given cell (in our case hyperbolic 24 cells) are equal.
If this is not the case, then we say the two horoballs are of different type.}

In order to compute volumes of horoball pieces, we use J\'anos Bolyai's classical formulas from the mid 19-th century:
\begin{enumerate}
\item
The hyperbolic length $L(x)$ of a horospheric arc that belongs to a chord segment of length $x$ is
\begin{equation}
\label{eq:horo_dist}
L(x)=2 \sinh{\left(\frac{x}{2}\right)}. \tag{2.9}
\end{equation}
\item The intrinsic geometry of a horosphere is Euclidean,
so the $(n-1)$-dimensional volume $\mathcal{A}$ of a polyhedron $A$ on the
surface of the horosphere can be calculated as in $\mathbb{E}^{n-1}$.
The volume of the horoball piece $\mathcal{H}(A)$ determined by $A$ and
the aggregate of axes
drawn from $A$ to the center of the horoball is
\begin{equation}
\label{eq:bolyai}
Vol(\mathcal{H}(A)) = \frac{1}{n-1}\mathcal{A}. \tag{2.10}
\end{equation}
\end{enumerate}
\section{On hyperbolic 24 cell}

An $n$-dimensional honeycomb $\mathcal{P}$, also referred to as a
solid tessellation or tiling, is an infinite collection of congruent
polyhedra (polytopes) that fit together face-to-face to fill the entire geometric
space $(at ~ present ~ \mathbb{H}^n~ (d \geqq 2))$ exactly once.
We take the cells to be congruent regular polyhedra. A honeycomb with cells congruent to a
given regular polyhedron $P$ exists if and only if the dihedral
angle of $P$ is a submultiple of $2\pi$ (in the hyperbolic plane
zero angles are also permissible). A complete classification of
honeycombs with bounded cells was first given by {{Schlegel}} in
$1883$. The classification was completed by including the polyhedra
with unbounded cells, namely the fully asymptotic ones by
{{Coxeter}} in 1954 \cite{C56}. {\it Such honeycombs (Coxeter tilings) exist only
for $d \le 5$} in hyperbolic $d$-space $\mathbb{H}^d$.

An alternative approach to describing honeycombs involves analysis
of their symmetry groups. If $\mathcal{P}$ is a Coxeter honeycomb,
then any rigid motion moving one cell into another maps the entire
honeycomb onto itself. The symmetry group of a honeycomb is denoted
by $Sym \mathcal{P}$. The characteristic simplex $\mathcal{F}$ of
any cell $P \in \mathcal{P}$ is a fundamental domain of the symmetry
group $Sym \mathcal{P}$ generated by reflections in its facets which
are $(d-1)$-dimensional hyperfaces.

The scheme of a regular polytope $P$ is a weighted graph (diagram)
characterizing $P \subset \mathbb{H}^d$ up to congruence. The nodes
of the scheme, numbered by $0,1,\dots,d$, correspond to the bounding
hyperplanes of $\mathcal{F}$. Two nodes are joined by an edge if the
corresponding hyperplanes are non-orthogonal. Let the set of weights
$(n_1,n_2,$ $n_3,\dots,n_{d-1})$ be the Schl\"afli symbol of $P$,
and $n_d$ be the weight describing the dihedral angle of $P$, such
that the dihedral angle is equal to $ \frac{2\pi}{n_d}$. In this
case $\mathcal{F}$ is the Coxeter simplex with the scheme:

\begin{figure}[ht]
\centering
\includegraphics[width=7cm]{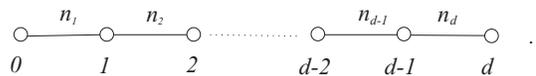}
\caption{Coxeter-Schl\"afli simplex scheme}
\end{figure}

The Schl$\ddot{a}$fli symbol of the honeycomb $\mathcal{P}$ is the
ordered set $(n_1,n_2,n_3,\dots,$ $n_{d-1},n_d)$ above.  A $(d+1)\times(d+1)$
symmetric matrix $(b^{ij})$ is constructed for each scheme in the
following manner: $b^{ii}=1$ and if $i \ne j\in \{0,1,2,\dots,d \}$
then $b^{ij}=-\cos{\frac{\pi}{n_{ij}}}$. For all angles between the
facets $i$,$j$ of $\mathcal{F}$ holds then $n_k= n_{k-1,k}$. Reversing the
numbering of the nodes of scheme $\mathcal{P}$ while keeping the
weights, leads to the scheme of the dual honeycomb $\mathcal{P}^*$
whose symmetry group coincides with $Sym \mathcal{P}$.

If $Sym\mathcal{P}$ denotes the symmetry group of a honeycomb then one tile $P$ of the Coxeter tiling
$\mathcal{P}_{n_1 n_2 \dots n_d}$ can be derived by the above symmetry group and its characteristic simplex $\mathcal{F}$:   
$$
P_{n_1 n_2 \dots n_d}=\Big\{\bigcup_{\gamma ~ \in ~ Sym \mathcal{P}_{n_1 n_2 \dots n_{d-1}}}
\gamma(\mathcal{F}_{n_1 n_2 \dots n_d}) \Big\}.
$$

Every $n$-dimensional totally asymptotic regular polytope $P$ has a hyperbolic ideal presentation obtained by normalising the
coordinates of its vertices so that they lie on the unit sphere $\mathbb{S}^{n-1}$ and by
interpreting $\mathbb{S}^{n-1}$ as the ideal boundary of $\overline{\mathbf{H}}^n$ in Beltrami-Cayley-Klein's ball model.
Therefore the ideal regular hyperbolic 24-cell $P_{24}$ can be derived from the Euclidean 24-cell as the convex hull
of the points
\begin{equation}
\begin{gathered}
A_1(1,\frac{1}{\sqrt{2}},\frac{1}{\sqrt{2}},0,0); \quad \quad \quad \quad A_{13}(1,-\frac{1}{\sqrt{2}},-\frac{1}{\sqrt{2}},0,0);\\
A_2(1,\frac{1}{\sqrt{2}},-\frac{1}{\sqrt{2}},0,0); \quad \quad \quad \quad A_{14}(1,-\frac{1}{\sqrt{2}},\frac{1}{\sqrt{2}},0,0);\\
A_3(1,\frac{1}{\sqrt{2}},0,\frac{1}{\sqrt{2}},0); \quad \quad \quad \quad A_{15}(1,-\frac{1}{\sqrt{2}},0,-\frac{1}{\sqrt{2}},0);\\
A_4(1,-\frac{1}{\sqrt{2}},0,\frac{1}{\sqrt{2}},0); \quad \quad \quad \quad A_{16}(1,\frac{1}{\sqrt{2}},0,-\frac{1}{\sqrt{2}},0);\\
A_5(1,\frac{1}{\sqrt{2}},0,0,\frac{1}{\sqrt{2}}); \quad \quad \quad \quad A_{17}(1,-\frac{1}{\sqrt{2}},0,0,-\frac{1}{\sqrt{2}});\\
A_6(1,-\frac{1}{\sqrt{2}},0,0,\frac{1}{\sqrt{2}}); \quad \quad \quad \quad A_{18}(1,\frac{1}{\sqrt{2}},0,0,-\frac{1}{\sqrt{2}});\\
A_7(1,0,\frac{1}{\sqrt{2}},\frac{1}{\sqrt{2}},0);  \quad \quad \quad \quad A_{19}(1,0,-\frac{1}{\sqrt{2}},-\frac{1}{\sqrt{2}},0);\\
A_8(1,0,-\frac{1}{\sqrt{2}},\frac{1}{\sqrt{2}},0); \quad \quad \quad \quad A_{20}(1,0,\frac{1}{\sqrt{2}},-\frac{1}{\sqrt{2}},0);\\
A_9(1,0,\frac{1}{\sqrt{2}},0,\frac{1}{\sqrt{2}}); \quad \quad \quad \quad A_{21}(1,0,-\frac{1}{\sqrt{2}},0,-\frac{1}{\sqrt{2}});\\
\end{gathered} \tag{3.1}
\end{equation}
\begin{equation}
\begin{gathered}
A_{10}(1,0,-\frac{1}{\sqrt{2}},0,\frac{1}{\sqrt{2}}); \quad \quad \quad \quad A_{22}(1,0,\frac{1}{\sqrt{2}},0,-\frac{1}{\sqrt{2}});\\
A_{11}(1,0,0,\frac{1}{\sqrt{2}},\frac{1}{\sqrt{2}}); \quad \quad \quad \quad A_{23}(1,0,0,-\frac{1}{\sqrt{2}},-\frac{1}{\sqrt{2}});\\
A_{12}(1,0,0,-\frac{1}{\sqrt{2}},\frac{1}{\sqrt{2}}); \quad \quad \quad \quad A_{24}(1,0,0,\frac{1}{\sqrt{2}},-\frac{1}{\sqrt{2}});
\end{gathered} \notag
\end{equation}
where the points (vertices) are described in a projective coordinate system given in Section 1.

The 24-cell is the unique regular four-dimensional polytope having cubical vertex figure because the vertex figure of the other five
regular four-dimensional polytopes are other Platonic solids, and therefore
their dihedral angles are not sub-multiples of $\pi$ thus only the regular 24-cell may be used as a
building block in order to construct cusped hyperbolic 4-manifolds.

\subsection{The structure of the hyperbolic 24-cell $P_{24}$}

$P_{24}$ is a tile of the 4-dimensional regular honeycomb $\mathcal{P}_{24}$ with Schl\"afli symbol $\{3,4,3,4 \}$.
It has 24 octahedral facets, 96 triangular faces, 96 edges and 24 cubical vertex figures.
A hyperbolic 24 cell contain $24 \cdot 48 = 1152$ characteristic simplex $\mathcal{F}_{24}$ and the volume of such a Coxeter simplex with
Schl\"afli symbol $\{3,4,3,4 \}$ is $Vol(\mathcal{F}_{24})=\frac{\pi^2}{864}$ (see \cite{JKRT})
therefore the volume of the hyperbolic 24-cell is $Vol(\mathcal{P}_{24})=\frac{4}{3}\pi^2$.

The vertices of $P_{24}$ are denoted by $A_i$ $(i\in \{1,2,\dots, 24\})$ and they coordinates are given in (3.1).

We introduce the notion of the $k$-neighbouring points $(k \in \{ 1,2,3,4 \})$ related to the vertices of $P_{24}$:
\begin{figure}[ht]
\centering
\includegraphics[width=10cm]{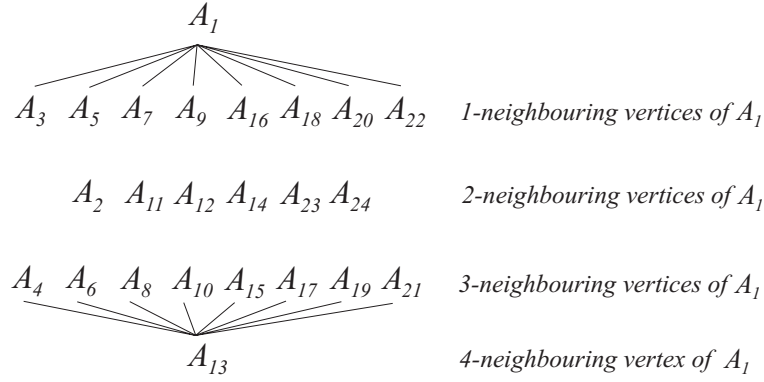}
\caption{The "neighborhood structure" of $P_{24}$}
\label{}
\end{figure}
\begin{defn}
\begin{enumerate}
\item The $1$-neighbouring vertices of $A_i$ $(i \in \{1,2,\dots, 24 \})$ among the vertices of $P_{24}$ are the vertices $A_j$
where $A_iA_j$ is an edge of $P_{24}$.
\item The $2$-neighbouring vertices of $A_i$ $(i \in \{1,2,\dots, 24 \})$ among the vertices of $P_{24}$ are $A_j$
where $A_iA_j$ is a diagonal of an octahedral facet of $P_{24}$.
\item The $4$-neighbouring vertex of $A_i$ $(i \in \{ 1,2,\dots, 24 \})$ among the vertices of $P_{24}$ is it opposite vertex $A_j$ regarding $P_{24}$.
\item The $3$-neighbouring vertices of $A_i$ $(i \in \{1,2,\dots, 24 \})$ among the vertices of $P_{24}$ are the vertices $A_j$
that are not $k$-neighbouring $(k=1,2,4)$ vertices of $A_i$.
\end{enumerate}
\end{defn}
The Fig.~2 shows the $k$-neighbouring vertices $(k=1,2,3,4)$ of $A_1$.
\begin{defn}
Two horoballs $B_i$ and $B_j$ (or horospheres $B_i^s$ and $B_j^s$ $(i,j \in \{1,2,\dots,24\},\ \ i \ne j)$ 
among the horoballs centered at the vertices $P_{24}$ 
are $k$-neighbouring $(k \in \{ 1,2,3,4 \})$ if their centres $A_i$ and $A_j$ are $k$-neighbouring vertices regarding $P_{24}$.
\end{defn}
We choose a characteristic
simplex (orthoscheme) of $P_{24}$ with vertices $T_0=A_1 \Big(1,\frac{1}{\sqrt{2}}, \frac{1}{\sqrt{2}},$ $0,0 \Big)$,
$T_1$, $T_2$, $T_3$ and $T_4=O$ where $T_4(1,0,0,0,0)$ is the centre of $P_{24}$
(coincides with the center of the model), $T_3 \Big(1,\frac{1}{2\sqrt{2}}, \frac{1}{2\sqrt{2}},\frac{1}{2\sqrt{2}},\frac{1}{2\sqrt{2}} \Big)$
is the centre of the facet-polyhedron $A_1 A_3 A_5 A_{7} A_{9} A_{11}$ (octahedron), the centre
of its regular face-polygon $A_1 A_3 A_{7}$ (regular triangle) is denoted by
$T_2 \Big(1,\frac{2}{3\sqrt{2}}, \frac{2}{3\sqrt{2}},\frac{2}{3\sqrt{2}},0 \Big)$ and $T_1
\Big(1,\frac{1}{\sqrt{2}}, \frac{1}{2\sqrt{2}},\frac{1}{2\sqrt{2}},0 \Big)$ is the centre of the edge $A_1A_3$ of this
face. Moreover, we denote by $T \Big(1,\frac{1}{2\sqrt{2}},$ $\frac{1}{2\sqrt{2}},\frac{1}{\sqrt{2}},0 \Big)$ the center of the edge $A_3A_7$.
This point is coincide with the orthogonal projection of $A_1$ onto its adjacent octahedral facet $A_3 A_4 A_7 A_8 A_{11} A_{24}$
(see Fig.~3).
\begin{figure}[ht]
\centering
\includegraphics[width=6cm]{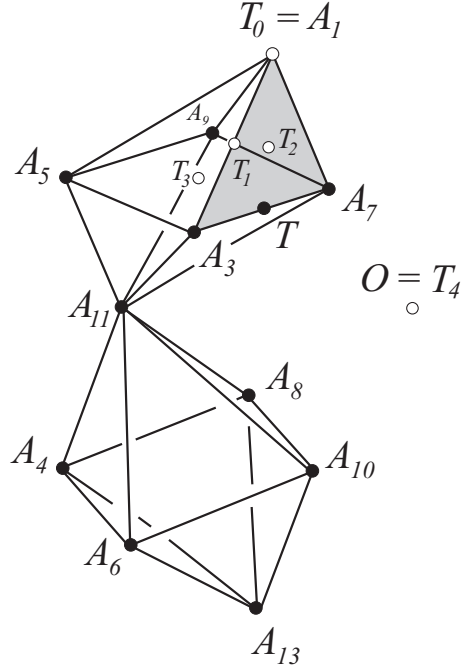}
\caption{A part of $P_{24}$}
\label{}
\end{figure}
\section{Horoball packings and polyhedral density function}
Similarly to the above section let $P_{24}$ be a tile of the 4-dimensional regular honeycomb $\mathcal{P}_{24}$ with Schl\"afli symbol $\{3,4,3,4 \}$. We study the horoball packings
$\mathcal{B}$ with horoballs centred at the infinite vertices of $\mathcal{P}_{24}$. 
The horospheres and horoball centred at the vertex $A_i$ are denoted by $B_i^s$ and $B_i$.
The density $\delta(\mathcal{B})$ of the horoball packing $\mathcal{B}$ relating to the above Coxeter tiling can be defined as the
extension of the local density related to the polytop $P_{24}$.
It is well known that for periodic ball or horoball packings the local density can be extended to the entire hyperbolic space.
\begin{defn}
We consider the polytop $P_{24}$ with vertices $A_i$ $(i=1, \dots, 24)$ in $4$-dimensional hyperbolic space $\overline{\mathbf{H}}^4$.
Centres of horoballs lie at vertices of $P_{24}$. We allow horoballs $(B_i, ~ i=0,1,2,\dots,24)$ of different types at the various vertices
and require to form a packing, moreover we assume that $$card \Big[ B_i \cap int \Big\{ \cup_{j=1}^{18} \mathcal{O}_{i_j} \Big\} \Big] = 0, $$
where the hyperplanes $\mathcal{O}_{i_j}$ $(j=1,\dots,18)$ do not contain the vertex $A_i$.
The generalized polyhedral density function for the above polytop and horoballs is defined as
\begin{equation}
\delta(\mathcal{B})=\frac{\sum_{i=0}^{24} Vol(B_i \cap P_{24})}{Vol(P_{24})}. \notag
\end{equation}
\end{defn}
The aim of this section is to determine the optimal packing arrangements $\mathcal{B}_{opt}$ and their
densities for the regular honeycomb $\cP_{24}$ in $\overline{\mathbf{H}}^4$. 
We vary the types of the horoballs so that they satisfy our constraints of non-overlap.
The packing density is obtained by the above definition.

We will use the consequences of the following Lemma (see \cite{Sz13}):

\begin{lemma}
Let $B_1$ and $B_2$ denote two horoballs with ideal centers $C_1$ and
$C_2$, respectively, in the $n$-dimensional hyperbolic space $(n \ge 2)$. Take $\tau_1$ and $\tau_2$ to be two congruent
$n$-dimensional convex piramid-like regions, with vertices $C_1$ and $C_2$. Assume that these horoballs
$B_1(x)$ and $B_2(x)$ are tangent at point $I(x)\in {C_1C_2}$ and
${C_1C_2}$ is a common edge of $\tau_1$ and
$\tau_2$. We define the point of contact $I(0)$ (the so-called ,,midpoint")such that the
following equality holds for the volumes of horoball sectors:
\begin{equation}
V(0):= 2 vol(B_1(0) \cap \tau_1) = 2 vol(B_2(0) \cap \tau_2). \notag
\label{szirmai-lemma}
\end{equation}
If $x$ denotes the hyperbolic distance between $I(0)$ and $I(x)$,
then the function
\begin{equation}
V(x):= vol(B_1(x) \cap \tau_1) + vol(B_2(x) \cap \tau_2)=\frac{V(0)}{2}(e^{(n-1)x}+e^{-(n-1)x}) \notag
\end{equation}
strictly increases as~$x\rightarrow\pm\infty$.
\end{lemma}
We consider the following four basic horoball configurations
$\mathcal{B}_i$, $(i=0,1,2,3,4)$:

\begin{enumerate}
\item All $24$ horoballs are of the same type and the adjacent horoballs
touch each other at the ,,midpoints" of each edge. This horoball arrangement is denoted by $\cB_0$.
\item We allow {\it horoballs of different types} and the opposite horoballs e.g. $B_1$ and $B_{13}$ touch their common 2-neighbouring horoballs $B_i$ 
$(i=2,11,12,14,23,$ $24)$ (see Fig.~2) at the centres of 
the corresponding octahedral facets e.g. the horoball $B_1$ touches the horoball $B_{11}$ at the facet center $T_3$ and 
$B_{13}$ tangent $B_{11}$ at the centre of octahedral facet $A_4 A_6  A_8 A_{10} A_{11} A_{13}$ (see Fig.~3). 
The other "smaller" horoballs are in the same type regarding $P_{24}$ and touch their 1-neighbouring "larger" horoballs e.g. 
the "larger" horoballs $B_1$ and $B_{11}$ touch the "smaller" horoballs $B_3,B_5,B_{7},B_{9}$. At this horoball arrangement let
the point $A_1A_3 \cap B_1^s$ be denoted by $C=I_1$ (see Fig.~4.a) ($B_i^s$ is the corresponding horosphere of horoball $B_i$.) 

This horoball arrangement is denoted by $\cB_1$. 
\item We set out from the $\cB_1$ ball configuration and we expand the horoballs $B_1$ and $B_{13}$ until they
comes into contact with their adjacent facets regarding $P_{24}$ while keeping their $1$ and $2$-neighbouring horoballs tangent to them.
At this configuration which is denoted by $\cB_2$ the horoballs are included on $3$ classes related to $P_{24}$.
The horoballs $B_1$ and $B_{13}$ are in the same type and they touch their corresponding $1$-neighbouring horoballs that form the second class.
The remaining $8$ horoballs are also in same type and are included on the $3$. type. 

For example the horoball $B_1$ touches its neighbouring facet at the point $T$ (see Fig.~3, and Fig.~4.b) and touches 
its $1$-neighbouring horoballs 
e.g. $B_3,B_5,B_{7},B_{9}$ and its $2$-neighbouring horoballs e.g. $B_{11}$. At this horoball arrangement let
the point $A_1A_{11} \cap B_1^s$ be denoted by $E=I_3$ (see Fig.~4.b). 
\item We set out also from the $\cB_1$ ball configuration and we expand the horoball $B_1$ 
until they comes into contact with their adjacent facets regarding $P_{24}$ while keeping their $1$ and 
$2$-neighbouring horoballs tangent to them. Moreover, we "blow up" the $3$-neighbouring horoballs of $B_1$ while their $1$-neigh\-bouring
horoballs touch them. At this configuration e.g. the horoball $B_1$ 
touches its neighbouring facet $A_3 A_4 A_{7} A_8 A_{11} A_{24}$ at the point $T$ (see Fig.~3, and Fig.~4.b) 
and touch its $1$-neighbouring horoballs e.g. $B_3,B_5,B_{7},B_{9}$ and its $2$-neighbouring horoballs e.g. $B_{11}$. 
Furthermore, the "expanded" horoballs e.g. $B_4, B_6, B_8, B_{10}$ touch the "shrunk" horoballs $B_{11}$ and $B_{13}$.  

This horoball arrangement is denoted by $\cB_3$.
\item Now we start from the configuration $\cB_0$ and we choose three arbitrary, mutually $3$-neighbouring horoballs and expand them
until they comes into contact with each other while keeping their $1$-neighbouring horoballs tangent to them.
We note here that this horoball configuration can be realized in $\bH^4$ (see the subsection 4.2.4). 
At this configuration which is denoted by $\cB_4$ the horoballs are included on $2$ classes related to $P_{24}$, e.g.
the horoballs $B_1$, $B_{10}$, $B_{17}$ are in same type touching each other and their "smaller" $1$-neighbouring horoballs that are also
in same type. 

\end{enumerate}
\subsection{Optimal horoball packings with horoballs in same type}
In this Section we consider the packings of horoballs where $Vol(B_i \cap P_{24})=Vol(B_j \cap P_{24})$ for all $i,j\in \{ 1,2, \dots, 24 \}$ 
thus the horoballs $B_i$ are in the same type regarding $P_{24}$. 

It is clear that in this case the maximal density can be achieved if the neighbouring horoballs touch each other at the centres of the edges 
of $P_{24}$ and the density of this densest packing $\mathcal{\cB}_0$ is equal to the maximal density of the horoball packings related to the 
Coxeter simplex tiling $\{3,4,3,4 \}$. 
For example in this case two horoballs $B_1$ and $B_3$ touch at the "midpoint" $T_1$ of edge $A_1A_3$ as projection of the polyhedron centre on 
it (see Fig.~3). 
These ball packings were investigated by the author in \cite{Sz07-1}: 
\begin{equation}
\begin{gathered}
V_0:=Vol(B_i \cap \mathcal{F}_{24})= \frac{1}{216}\sqrt{2}\sinh{\Big( \frac{1}{2} \mathrm{arcosh}{\Big(\frac{11}{8}\Big)}\Big)}\approx 0.00694, \\
Vol(\mathcal{F}_{24})=\frac{\pi^2}{864}, \ \
\delta(\mathcal{B}_0)=\frac{Vol(B_i \cap \mathcal{F}_{24})}{Vol(\mathcal{F}_{24})}\approx 0.60793.
\end{gathered} \tag{4.1}
\end{equation}
\subsection{Optimal horoball packings with horoballs in different types}
The type of a horoball is allowed to expand until either the horoball
comes into contact with other horoballs or with a adjacent facet of the
honeycomb. These conditions are satisfactory to ensure
that the balls form a non-overlapping horoball arrangement, as such the collection of all
horoballs is a well defined packing in $\mathbb{H}^4$.
\subsubsection{Horoball packings $\cB_0^1$ and their densities between the horoball arrangements $\cB_0$ and $\cB_1$}
We set out from the $\cB_0$ ball configuration (see above Section) and consider two $1$-neighbouring horoballs e.g. $B_1$ and $B_{3}$ from it.
Let $I_0=I(0)=T_1$ be their point of tangency on side $A_1A_3$ (see Fig.~3 and Fig.~4.a). 
Moreover, consider the point $I(x)$ on the segment $A_1A_3$ 
where the modified horoballs $B_i(x), \ (i=1,3)$ are tangent to each other
and $x$ is the hyperbolic distance between $I(0)$ and $I(x)$ (the value of $x$ can also be negative if $I(x)$ is on the segment $T_1A_1$). 

We blow up the horoballs $B_1(0)$ and $B_{11}(0)$ (and also the horoballs $B_2$, $B_{12}$, $B_{14}$, $B_{23}$, $B_{24}$ and $B_{13}$
to achieve the $\cB_1$ horoball configuration) until they come 
into contact with each other at the centre $T_3$ of octahedral facet $A_1 A_3 A_{5} A_7 A_{9} A_{11}$. At this situation 
(see Fig.~3) the horoball centered at $A_1$ is denoted, by $B_1(\rho_1)$ where $\rho_1$ is the hyperbolic distance between $I_0$ and $I_1$ (see Fig.~4.a).

The foot-point of the perpendicular from $T_3$ onto the staight line $A_1A_3$ is $I_0=T_1$ which is the common point of the horoballs
$B_1(0) \in \cB_0$ and $B_{3}(0) \in \cB_0$ centered at $A_1$ and $A_{3}$, respectively. The hyperbolic distance $s_1=T_1T_3$ 
between the points $T_1[\bt_1]$ and $T_3[\bt_3]$ can be computed by the formula (2.5) (see Fig.~4.a):
\begin{figure}[ht]
\centering
\includegraphics[width=11cm]{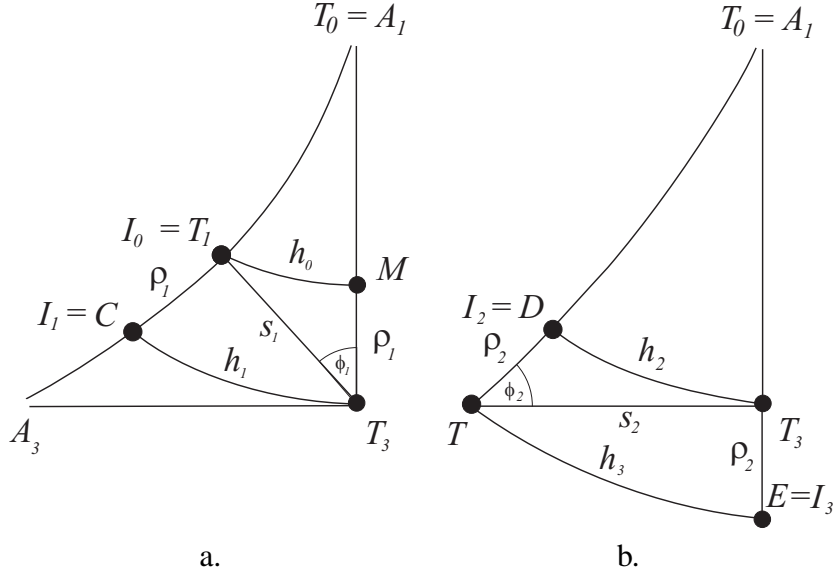}

a. \hspace{5cm} b.
\caption{Computations of hyperbolic distances $\rho_1=I_0I_1$ and $\rho_2=I_2T=I_3T_3$}
\label{}
\end{figure}
The parallel distance of the angle $\phi_1=T_1T_3A_1\angle$ is $s_1$ therefore we obtain by the classical formula of J.~Bolyai and by formula (2.5) the
following equation (see Fig.~4.a).
\begin{equation}
\frac{1}{\sin{(\phi_1})}=\cosh{s_1}=\sqrt{2}. \tag{4.2}
\end{equation}

We consider two horocycles $\mathcal{H}_0$ and $\mathcal{H}_1$ through the points $I_0$ and $I_1$ with center $A_1$ in the plane $A_1A_3T_3$ 
and the point $\mathcal{H}_1 \cap A_1T_3$ is denoted by $M$.
The horocyclic distances between points $I_0$, $M$ and $I_1$, $T_3$ are denoted by $h_0$ and $h_1$. 
By means of formula of J.~Bolyai and of (4.2), we have
\begin{equation}
\begin{gathered}
\frac{h_1}{h_0}=e^{\rho_1}=\frac{1}{\sin(\phi_1)} \ \Rightarrow \rho_1=\log(\sqrt{2}) \approx 0.34657.
\end{gathered} \tag{4.3}
\end{equation}
We extend the above modifications and denotations for all horoballs of packings between horoball arrangements $\cB_0$ and $\cB_1$ i.e. the horoballs are denoted by $B_i(x)$
$(i \in [0,\rho_1])$. If $x=0$ then we get the $\cB_0$ horoball packing and if $x=\rho_1$ then the $\cB_1$ one.

We obtain using the results of the former computations and of Lemma 4.2  the next 
\begin{lemma} 
The density of packings $\cB_0^1$ (see Fig.~5.a) between the 
main horoball arrangements $\cB_0$ and $\cB_1$ can be computed by the formula 
\begin{equation}
\begin{gathered}
\delta(\mathcal{B}_0^1(x))=\frac{\sum_{i=0}^{24} Vol(B_i(x) \cap P_{24})}{Vol(P_{24})}= 
\frac{384 \cdot V_0~(e^{3x} + 2 \cdot e^{-3x})}{\frac{4}{3}\pi^2},~ ~ x \in [0,\rho_1] 
\end{gathered} \notag
\end{equation}
and the maxima of function $\delta(\mathcal{B}_0^1(x))$ (see Fig.~5.a) are realized at $x=\rho_1 ~\approx 0.34657$ where the horoball packing density
is $\delta(\mathcal{B}_0^1(\rho_1))\approx 0.71645$. 
\end{lemma}
\begin{figure}[ht]
\centering
\includegraphics[width=6cm]{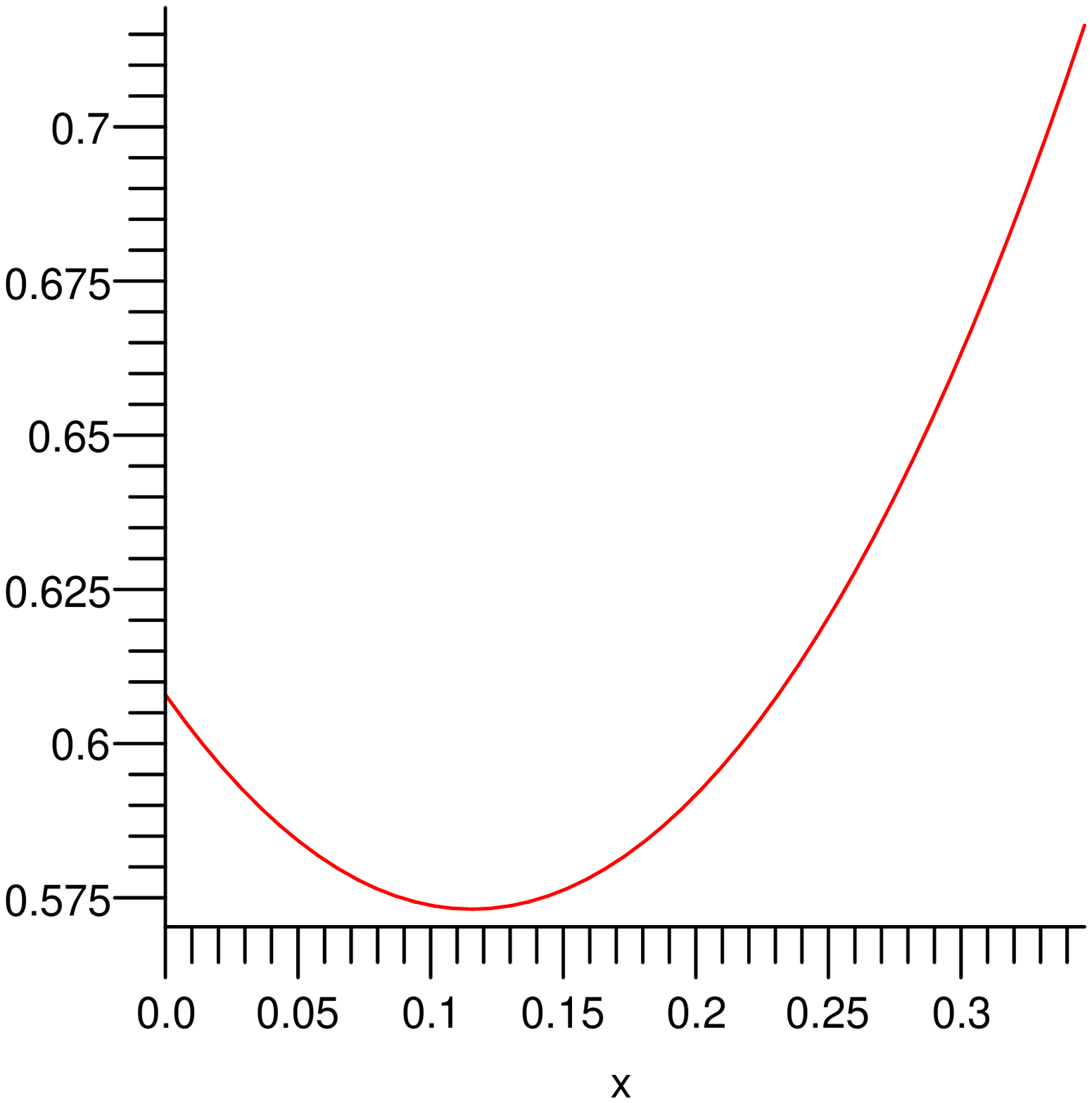} \includegraphics[width=6cm]{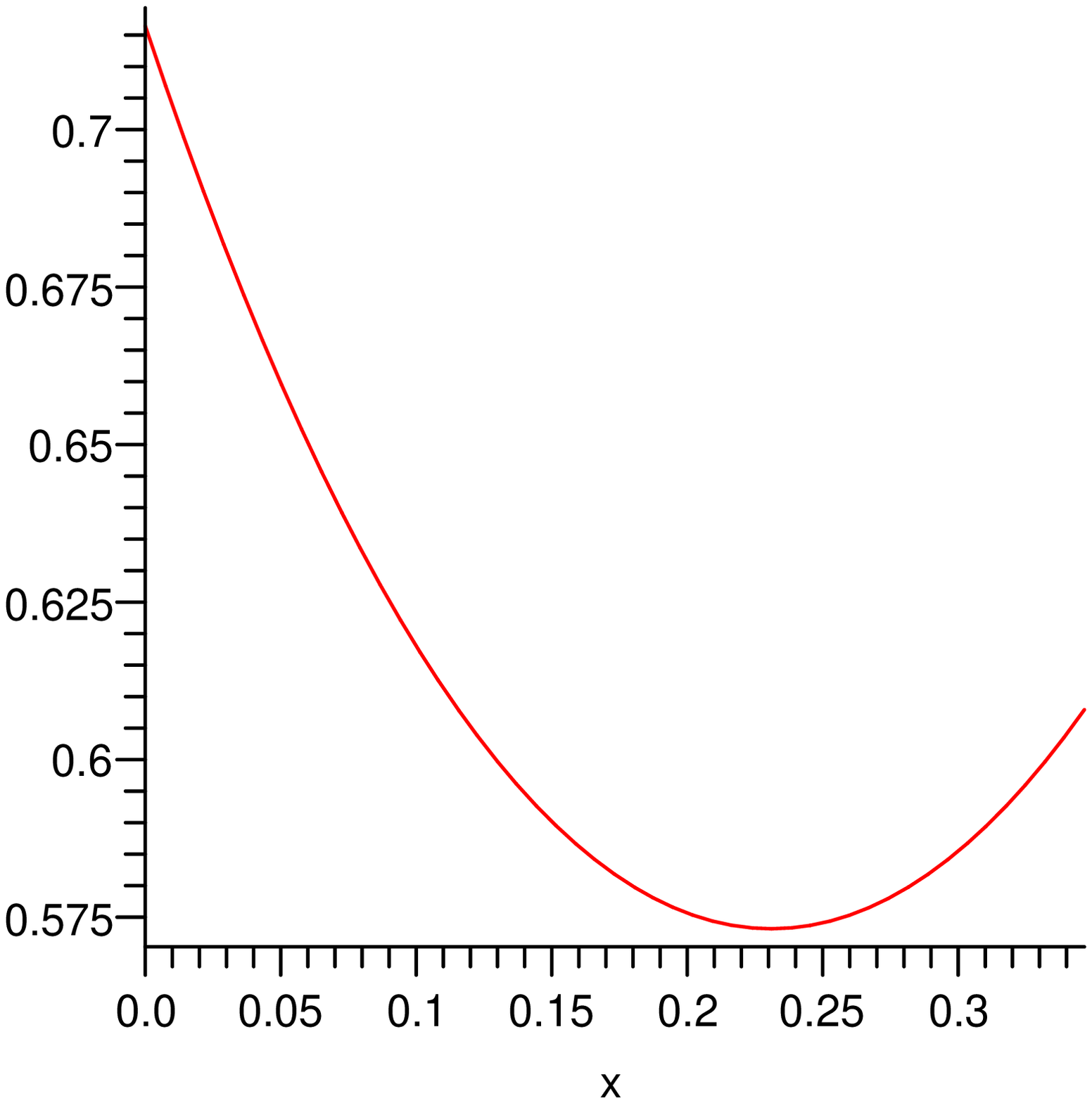}

a. \hspace{5cm} b.
\caption{The graphs of functions $\delta(\mathcal{B}_0^1(x))$ and $\delta(\mathcal{B}_1^2(x))$ where $x \in [0,\rho_1]$.}
\label{}
\end{figure}
\begin{rmrk}
We note, here that the above optimal density $\delta(\mathcal{B}_0^1(\rho_1))\approx 0.71645$ is equal to the density 
of known densest ball and horoball packings in $\mathbb{H}^4$.
\end{rmrk}
\subsubsection{Horoball packings $\cB_1^2$ and their densities between the horoball arrangements $\cB_1$ and $\cB_2$}
We start our investigation from the $\cB_1$ ball configuration. Here 
e.g. the horoballs $B_1$ and $B_{3}$ touch each other at the point $I_1$ (see Fig.~4.a) and $B_1$ touch $B_{11}$ at the point $T_3$ 
(see Fig.~3 and Fig.~4.b) etc. Furthermore, in $\cB_1$ the common point of horosphere $B_1^s$ with the line segment $A_1T$ is denoted by 
$I_2=I^*(0)=D$ (see Fig.~4.b). We consider the point $I^*(x)$ on the segment $A_1T$ 
where a modified horosphere $B_1^s(x)$ intersects the line segment $A_1T$ 
and $x$ is the hyperbolic distance between $I^*(0)$ and $I^*(x)$ (the value of $x$ can also be negative if $I^*(x)$ is on the segment $A_1I^*(0)$). 
Corresponding to the above notions we introduce the notations $B_{13}(0)$ and $B_{13}(x)$. 

We blow up the horoballs $B_1(0)$ and $B_{13}(0)$ while keeping 
their $1$-neighbouring horoballs tangent to them until they comes 
into contact with their adjacent facets of $P_{24}$ e.g. upto the horoball $B_1(x)$ touches the octahedral facet $A_3 A_4 A_{7} A_8 A_{11} A_{24}$. 
At this arrangement relating to Fig.~4.b the horoball centered at $A_1$ is denoted, by $B_1(\rho_2)$ 
where $\rho_2$ is the hyperbolic distance between $I^*(0)$ and $T$.

The foot-point of the perpendicular from $T$ onto the staight line $A_1A_{11}$ is $T_3$. 
The hyperbolic distance $s_2=T T_3$ 
between the point $T[\bt]$ and $T_3[\bt_1]$ can be computed by the formula (2.5) (see Fig.~4.b):
The parallel distance of the angle $\phi_2=A_1TT_3\angle$ is $s_2$ therefore we obtain by the classical formula 
of J.~Bolyai and by formula (2.5) the
following equation (see Fig.~4.b):
\begin{equation}
\frac{1}{\sin{(\phi_2})}=\cosh{s_2}=\sqrt{2}. \tag{4.4}
\end{equation}

We consider two horocycles $\mathcal{H}_2$ and $\mathcal{H}_3$ through the points $I_2$ and $T$ with center $A_1$ in the plane $A_1TT_3$ 
and the point $\mathcal{H}_3 \cap A_1T_3$ is denoted by $E=I_3$.
The horocyclic distances between points $I_2$, $T_3$ and $T$, $E$ are denoted by $h_2$ and $h_3$. 
Similarly to (4.3) we obtain that $\rho_2=\log(\sqrt{2}) \approx 0.34657$.

We extend the above modifications and denotations for all horoballs of packings between horoball arrangements $\cB_1$ and $\cB_2$ i.e. the horoballs are denoted by $B_i(x)$
$(i \in [0,\rho_2])$. If $x=0$ then we get the $\cB_1$ horoball packing and if $x=\rho_2$ then the $\cB_2$ one.

We obtain using the results of the former computations and of Lemma 4.2  the next 
\begin{lemma} 
The density of packings $\cB_1^2$ (see Fig.~5.b) between the 
main horoball arrangements $\cB_1$ and $\cB_2$ can be computed by the formula 
\begin{equation}
\begin{gathered}
\delta(\mathcal{B}_1^2(x))=\frac{\sum_{i=0}^{24} Vol(B_i(x) \cap P_{24})}{Vol(P_{24})}= \\
= \frac{48 \cdot V_0~(2~ e^{3(\rho_1+x)} + 6\cdot e^{-3(-\rho_1+x)}+16 \cdot e^{-3(\rho_1+x)})}{\frac{4}{3}\pi^2},~ ~ x \in [0,\rho_2] 
\end{gathered} \notag
\end{equation}
and the maxima of function $\delta(\mathcal{B}_1^2(x))$ (see Fig.~5.b) are realized at $x=0$ i.e. at the $\cB_1$ ball packing (see Lemma 4.3).
\end{lemma}
\begin{rmrk}
The density $\delta(\mathcal{B}_1^2(\rho_2))$ is equal to the maximal density of packing with horoballs in same types: 
$\delta(\mathcal{B}_1^2(\rho_2))=\delta(\mathcal{B}_0) \approx 0.60793$.
\end{rmrk}
\subsubsection{Horoball packings $\cB_1^3$ and their densities between the horoball arrangements $\cB_1$ and $\cB_3$}
Similarly to the above subsection we set out from the $\cB_1$ ball configuration and we will use the notations of subsection 4.2.2. 
Now, we expand the horoball $B_1(0)$  
until they comes into contact with their adjacent facets regarding $P_{24}$ while keeping their $1$ and 
$2$-neighbouring horoballs tangent to them. Moreover, we "blow up" the $3$-neighbouring horoballs of $B_1(0)$ while their $1$-neighbouring
horoballs touch them. At this procedure this horoball is denoted by $B_1(x)$. If we achieved the endpoint of this extension then
e.g. the horoball $B_1(\rho_2)$ touches its neighbouring facet $A_3 A_4 A_{7} A_8 A_{11} A_{24}$ at the point $T$ (see Fig.~3, and Fig.~4) 
and touch its $1$-neighbouring horoballs e.g. $B_3,B_5,B_{7},B_{9}$ and its $2$-neighbouring horoballs e.g. $B_{11}$. 
Furthermore, the "expanded" horoballs e.g. $B_4, B_6, B_{8}, B_{10}$ touch the "shrunk" horoballs $B_{11}$ and $B_{13}$.

We extend the above modifications and notations for all horoballs of packings between horoball arrangements 
$\cB_1$ and $\cB_3$ i.e. the horoballs are denoted by $B_i(x)$
$(i \in [0,\rho_2])$. If $x=0$ then we get the $\cB_1$ horoball packing and if $x=\rho_2$ then the $\cB_3$ one.
Finally, we obtain the next
\begin{lemma} 
The density of packings $\cB_1^3$ between the 
main horoball arrangements $\cB_1$ and $\cB_3$ can be computed by the formula 
\begin{equation}
\begin{gathered}
\delta(\mathcal{B}_1^3(x))=\frac{\sum_{i=0}^{24} Vol(B_i(x) \cap P_{24})}{Vol(P_{24})}= \\
= \frac{48 \cdot V_0~(e^{3(\rho_1+x)} + 7 \cdot e^{-3(-\rho_1+x)}+ 8 \cdot e^{-3(\rho_1+x)}+8 \cdot e^{-3(\rho_1-x)})}{\frac{4}{3}\pi^2},~ x \in [0,\rho_2] 
\end{gathered} \notag
\end{equation}
and the maxima of function $\delta(\mathcal{B}_1^3(x))$ are realized at $x=0$ i.e. at the $\cB_1$ ball packing (see Lemma 4.3).
\end{lemma}
\begin{rmrk}
The function $\delta(\mathcal{B}_1^3(x))$ is the same with $\delta(\mathcal{B}_1^2(x))$ $(x \in [0,\rho_2])$ (see Fig.~5.b).
\end{rmrk}
\subsubsection{Horoball packings $\cB_0^4$ and their densities between the horoball arrangements $\cB_0$ and $\cB_4$}
Here we consider the horoball configuration $\cB_0$ and we choose three arbitrary, mutually $3$-neighbouring horoballs e.g. $B_1$,
$B_{10}$ and $B_{17}$ and let $I_6=I_*(0)$ be the point of intersection of horosphere $B_1^s(0)$ with the segment $T A_1$.
Moreover, consider the point $I_*(x)$ on the segment $I_6T$ 
where the expanded horosphere $B_1^s(x)$ intersects the segment $I_6T$ and $x$ is the 
hyperbolic distance between $I_*(0)$ and $I_*(x)$ (see Fig.~6.a). We have seen in former subsections 
that the hyperbolic distance between $I_0$ and $T$ is $2\rho_1=2\rho_2$ (see Fig.~5a and Fig.~5b). 
We consider a horocycles $\mathcal{H}_5$ through the point $T$ 
with center $A_1$ in the plane $A_1A_{10}T$ and the point $\mathcal{H}_5 \cap A_1A_{10}$ is denoted by $K=I_5$.

The foot-point of the perpendicular from $T$ onto the staight line $A_1A_{10}$ is called by $Q$ whose coordinates are 
$Q\Big(1,\frac{5}{7\sqrt{2}}, \frac{3}{7\sqrt{2}},0, \frac{2}{7\sqrt{2}}\Big)$. 

We obtain by in the subsections 4.2.1 and 4.2.2 described method that the hyperbolic distance $\rho_3$ of the points $Q$ and $K$ is
$\rho_3=\log{\frac{10}{3}}$. 

The centre ("midpoint") of segment $A_1A_{10}$ is denoted by $H$ (see Fig.~6.a) 
(in our model this is Euclidean midpoint of segment $A_1A_{10}$, 
as well) whose distance $\rho_4$ to $Q$ can be computed by the formula (2.5): $\rho_4={\text{arccosh}}\Big(\frac{7\sqrt{2}}{4\sqrt{5}}\Big)$. 
The point $H$ lie on the line segment $QK$ because $0.60199  \approx \rho_3 > \rho_4 \approx 0.45815$.  

Finally, we obtain using the results of the former computations and of Lemma 4.2  the next 
\begin{lemma} 
The density of packings $\cB_0^4$ (see Fig.~5.b) between the 
main horoball arrangements $\cB_0$ and $\cB_4$ can be computed by the formula 
\begin{equation}
\begin{gathered}
\delta(\mathcal{B}_0^4(x))=\frac{\sum_{i=0}^{24} Vol(B_i(x) \cap P_{24})}{Vol(P_{24})}= \\
= \frac{48 \cdot V_0~(3\cdot e^{3x} + 21\cdot e^{-3x})}{\frac{4}{3}\pi^2},~ ~ x \in [0,2\rho_1+\rho_4-\rho_3 \approx 0.54931] 
\end{gathered} \notag
\end{equation}
and the maxima of function $\delta(\mathcal{B}_0^4(x))$ (see Fig.~6.b) are realized at $x=0$ where the horoball packing density
is $\delta(\mathcal{B}_0^4(0))\approx 0.60793$. 
\end{lemma}
\begin{figure}[ht]
\centering
\includegraphics[width=4cm]{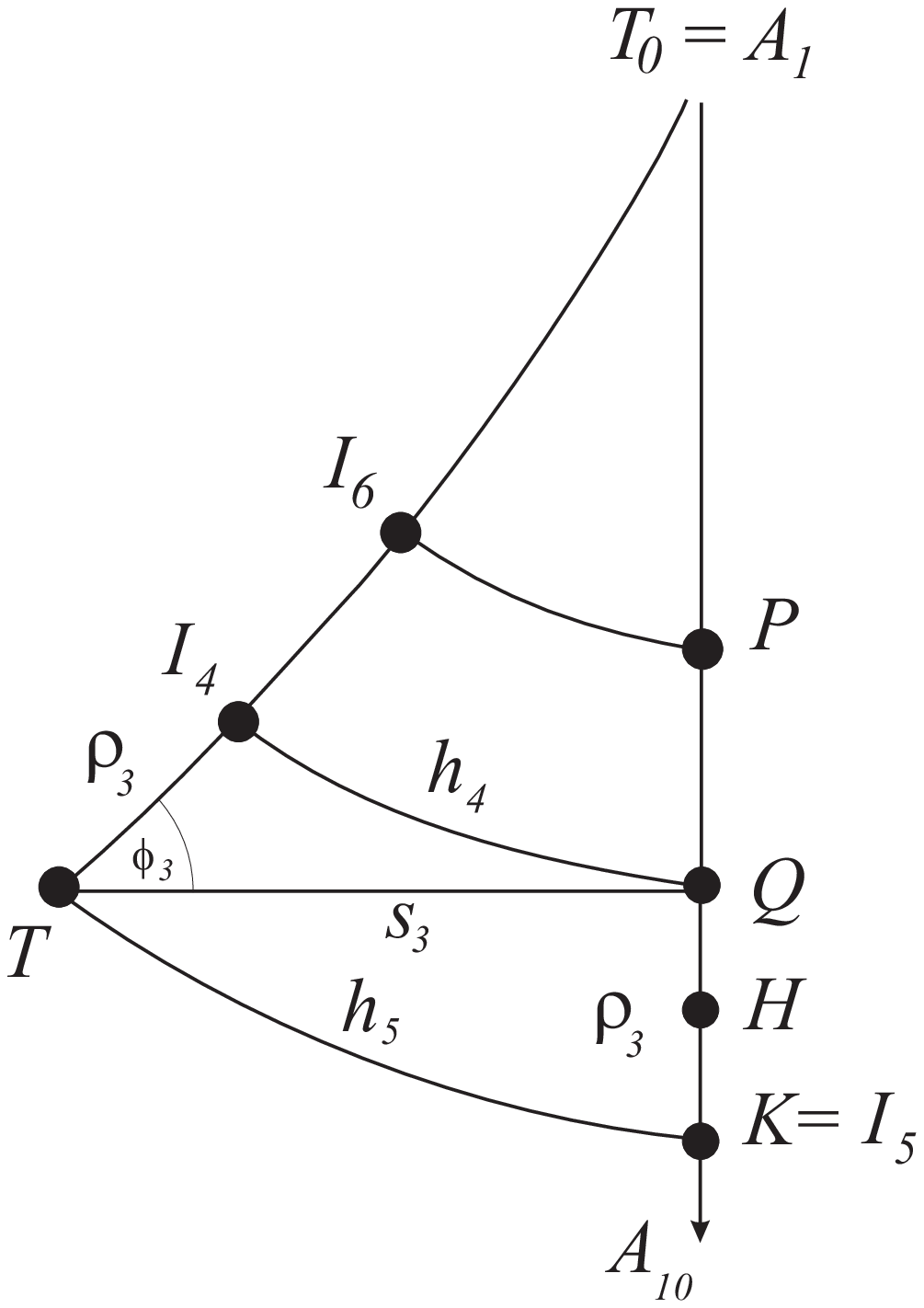} \includegraphics[width=6cm]{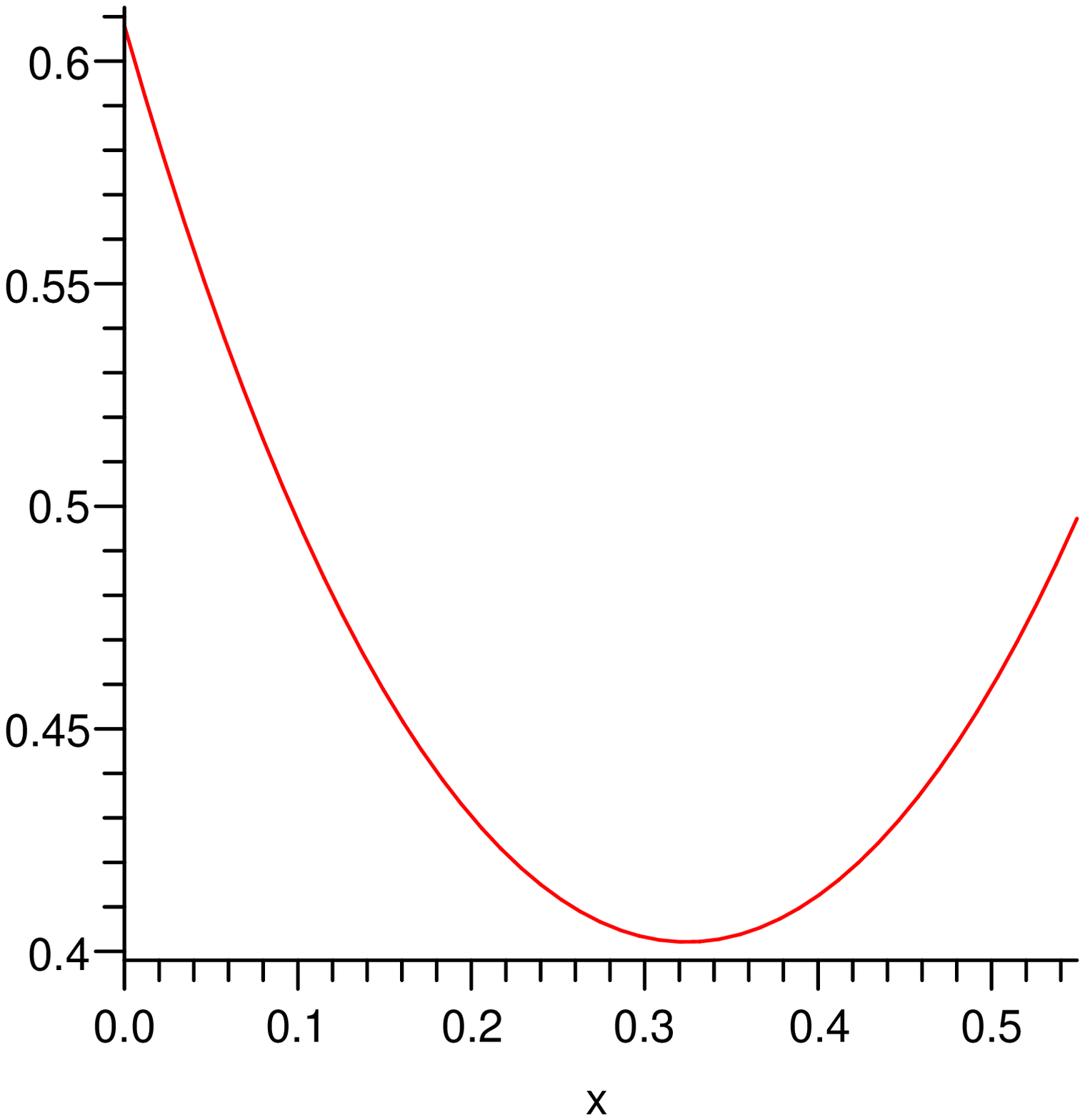}

a. \hspace{5cm} b.
\caption{a.~The computation of hyperbolic distance $\rho_3=I_4T$ and $\rho_4=QH$.~ b.~ 
The graph of function $\delta(\cB_0^4(x))$ $x \in [0, 2 \rho_1+ \rho_3- \rho_4 ]$.}
\label{}
\end{figure}
\begin{rmrk}
The density $\delta(\mathcal{B}_0^4(0))$ is equal to the maximal density of packings with horoballs in same types: 
$\delta(\mathcal{B}_0^4(\rho_2))=\delta(\mathcal{B}_0) \approx 0.60793$.
\end{rmrk}
\subsection{Optimal horoball packings to hyperbolic $24$-cell}
The main result of this paper is summarized in the following
\begin{theorem}
\label{mainresult}
The horoball arrangement $\cB_1$ (see 4.2.1) provide the maximal horoball packing density 
related to the hyperbolic tiling $\mathcal{P}_{24}$ with Schl\"afli symbol $\{3,4,3,4 \}$ and its 
density is $\delta_{opt}(\mathcal{B}) \approx 0.71645$ 
if horoballs of different types are allowed at each asymptotic vertex of the tiling. 
\end{theorem}
\begin{rmrk}
The optimal horoball packing described and determined in this paper is a new horoball configuration which provide the known 
maximal density of realizable packings of 
the entire hyperbolic space $\mathbb{H}^4$. 
\end{rmrk}
{\bf Proof}

It is well known that a packing is optimal, then it is locally stable i.e. each ball is fixed by the other ones so that no ball of packing 
can be moved alone without overlapping another ball of the given ball packing.

The packings of horoballs can be easily classified by the type of "maximally large" horoball regarding the horoball packing to $\mathcal{P}_{24}$.
If we fix the "maximally large" horoball related to the above tiling then all possible horoball packing can be modified to achieve 
one of the above horoball configurations $\cB_i^j(x)$ $(i,j\in {\{0,1,2,3,4\}}, ~ i<j)$ without decrease of the packing density.

A horoball $B_i(x)$ is "maximally large" if $Vol(B_i(x) \cap P_{24})$ $(i \in 1 \dots 24)$ is maximal. Here the maximal volume is denoted by 
$Vol(B_i^{max})$. 
\begin{enumerate}
\item If $\frac{1}{48}Vol(B_i^{max}) \le V_0$ then the maximal density can be computed by Sect. 4.1 where the maximal density is $\delta(\mathcal{B}_0) \approx 0.60793$. 
\item If $V_0 < \frac{1}{48} \cdot Vol(B_i^{max}) \le V_0 \cdot e^{3\rho_1}$ then the optimal density can be computed 
by Sections 4.2.1, here the optimal 
density is $\delta(\mathcal{B}_0^1(\rho_1))\approx 0.71645$.
\item If $V_0 \cdot e^{3\rho_1} < \frac{1}{48} \cdot Vol(B_i^{max}) \le V_0 \cdot e^{6\rho_1}$ then the densities can be computed by Sections 4.2.2, 
4.2.3 and 4.2.4 where the maximal density is $\delta(\mathcal{B}_0) \approx 0.60793$. 
\end{enumerate}

The volume of the "largest horoball" $Vol(B_i^{max}) \le V_0 \cdot e^{6\rho_1}$ therefore we proved the above Theorem. $\square$ 

\vspace{3mm}

The above results also show, that the discussion of the densest horoball packings and coverings in 
the $n$-dimensionalen hyperbolic space with horoballs
of different types has not been settled yet. Similarly to these, the problems of the densest hypersphere
(or hyperball) packings and coverings are open, as well.

Optimal sphere packings in other homogeneous Thurston geometries
form also a class of open mathematical problems (see 
\cite{Sz07-2}, \cite{Sz10-1}, \cite{Sz11-1}, \cite{Sz10-3}, \cite{Sz12-1}, \cite{Sz13-2}, \cite{MSz14}, 
\cite{MSzV14}, \cite{MSzV15}). Detailed studies are the objective of
ongoing research.



\end{document}